\documentclass[nologo,notitlepageheader]{dakota-article}

\usepackage{lineno,hyperref}
\usepackage{epsfig}

\usepackage{graphicx,amsfonts,amsmath,amssymb,mathrsfs,bbm}

\newcommand{\real}{\mathbbm{R}}
\newcommand{\complex}{\mathbbm{C}}
\newcommand{\nat}{\mathbbm{N}}
\newcommand{\ltwo}{\mathscr{L}^2(\Pi,\rho)}
\newcommand{\ltwotime}{\mathscr{L}^2[0,\infty)}

\title{Time and Frequency Domain Methods for Basis\- Selection in Random Linear Dynamical Systems}
\def\SNL{Center for Computing Research, Sandia National Laboratories,
  P.O. Box 5800, MS 1318, Albuquerque, NM, 87185-1320, United States.}
\def\GRE{Institute of Mathematics and Computer Science, University of Greifs\-wald,  Walther-Rathenau-Stra{\ss}e~47, D-17489 Greifswald, Germany.}
\shorttitle{Time and Frequency Domain Methods for Basis Selection in Random Linear Dynamical Systems}
\shortauthor{J.D. Jakeman \& R.~Pulch}
\author{John D. Jakeman\thanks{\SNL}, Roland Pulch\thanks{\GRE}}
\corauthor{Roland Pulch}
\coremail{roland.pulch@uni-greifswald.de}
\date{\today}
\funding{See acknowledgements}
 
\graphicspath{{../}}%delete
\begin{document}
\pagestyle{dakotaheader}
\maketitle

% Structure definitions
%\newtheorem{theorem}{Theorem}
%\newtheorem{lemma}{Lemma}
%\newtheorem{definition}{Definition}

%%%%%%%%%%%%%%%%%%%%%%%%%%%%%%%%%%%%%%%%%%%%%%%%%%%%%%%%%%%%%%%%%%%%%%%%%%%%%
%%%                         Abstract                                      %%%
%%%%%%%%%%%%%%%%%%%%%%%%%%%%%%%%%%%%%%%%%%%%%%%%%%%%%%%%%%%%%%%%%%%%%%%%%%%%%

\begin{abstract}
  Polynomial chaos methods have been extensively used to analyze
  systems in uncertainty quantification.
  Furthermore, several approaches exist to determine a 
  low-dimensional approximation (or sparse approximation)
  for some quantity of interest in a model,
  where just a few orthogonal basis polynomials are required.
  We consider linear dynamical systems consisting of ordinary differential
  equations with random variables.
  The aim of this paper is to explore methods for producing
  low-dimensional approximations of the quantity of interest
  further. 
  We investigate two numerical techniques to compute a low-dimensional
  representation, which both fit the approximation to a set of samples
  in the time domain.
  On the one hand, 
  a frequency domain analysis of a stochastic Galerkin
  system yields the selection of the basis polynomials. 
  It follows a linear least squares problem.
  On the other hand,
  a sparse minimization yields the choice of the basis polynomials
  by information from the time domain only. 
  An orthogonal matching pursuit produces an approximate solution of
  the minimization problem.
  We compare the two approaches using a test example from a
  mechanical application.
  \end{abstract}


\maketitle

%%%%%%%%%%%%%%%%%%%%%%%%%%%%%%%%%%%%%%%%%%%%%%%%%%%%%%%%%%%%%%%%%%%%%%%%%%%%%
%%%                       Introduction                                    %%%
%%%%%%%%%%%%%%%%%%%%%%%%%%%%%%%%%%%%%%%%%%%%%%%%%%%%%%%%%%%%%%%%%%%%%%%%%%%%%

\section{Introduction}
We consider linear dynamical systems in the 
form of ordinary differential equations (ODEs),
which include physical parameters.
A quantity of interest (QoI) is defined as an output of the system.
Uncertainties may be present in the parameters.
In uncertainty quantification, a common approach is to
interpret the parameters as random variables, 
see~\cite{sullivan,xiu-book}.

The state variables as well as the QoI can be
expanded into a series with given orthogonal basis polynomials
depending on the random variables and a priori unknown time-dependent
coefficient functions.
This spectral approach is called a (generalized) polynomial chaos expansion,
see~\cite{augustinetal,ernst,sullivan,xiu-book}.
Several types of numerical methods exist to compute an approximation of
the coefficient functions.
On the one hand, 
the stochastic Galerkin method projects the random linear dynamical system
to a larger deterministic linear system of ODEs, whose solution
yields the approximation, see~\cite{pulch11,pulch14}.
On the other hand,
the approximation can be fitted to random samples in a
least squares regression, see~\cite{hampton-doostan}.

The number of basis polynomials up to a given total degree becomes huge
in the case of a large number of random variables.
Our task is to identify a low-dimen\-sional approximation of the
random QoI, which is sometimes called a sparse representation
in the literature.
Therein, just a small subset of basis polynomials is required for
a sufficiently accurate approximation.
Numerical methods for this problem were based on 
least angle regression~\cite{blatman}, 
sparse grid quadrature~\cite{conrad-marzouk}, 
compressed sensing~\cite{doostan,Yan_GX_IJUQ_2012}, 
and reduced basis techniques \cite{nair,sachdeva}, for example.
Dimension-adaptive ANOVA decompositions~\cite{Prasad_R_IEEETMTT_2017,Yang_CLK_JCP_2012,Zhang_YOKD_IEEECDICS_2015}, and low-rank tensor approximations~\cite{Oseledets_2013,Zhang_WD_IEEETCPMT_2017} can also be used to address the curse-of-dimensionality.

In this paper, the aim is to further explore numerical methods for
low-dimensional representations. 
We investigate and compare two techniques to construct a
low-dimensional approximation with $q$~basis polynomials,
where~$q$ is a given integer number. 
Both approaches fit their approximations to a set of samples from the QoI
in the time domain. 
First, a frequency domain analysis is performed for the transfer function
of the stochastic Galerkin system,
which was derived in the previous work~\cite{pulch2018}. 
The minimization of an error bound identifies a subset of
$q$~basis polynomials.
We apply this subset and identify the accompanying
time-dependent coefficient functions by a least squares problem,
which was not considered in~\cite{pulch2018}. 
In addition, we examine an $\ell_0$-minimization under a constraint
purely in the time domain,
which also yields a subset of $q$~basis polynomials.
Orthogonal matching pursuit (OMP)~\cite{Pati1993}
generates a numerical approximation to the minimization problem.

We apply both techniques to a test example, which models a
mass-spring-damper system with random parameters.
A comparison of the errors for the low-dimen\-sional approximations
is presented.
Moreover, we examine the conformance of both approaches, i.e.,
if the techniques predominantly identify the same basis polynomials.

%%%%%%%%%%%%%%%%%%%%%%%%%%%%%%%%%%%%%%%%%%%%%%%%%%%%%%%%%%%%%%%%%%%%%%%%%%%%%
%%%                      Problem Definition                               %%%
%%%%%%%%%%%%%%%%%%%%%%%%%%%%%%%%%%%%%%%%%%%%%%%%%%%%%%%%%%%%%%%%%%%%%%%%%%%%%

\section{Problem Definition}
\label{sec:problem}
The problem of the basis selection is specified in this section.
We denote the sets of real numbers and complex numbers by
$\real$ and $\complex$, respectively. 

%%%%%%%%%%%%%%%%%%%%%%%%%%%%%%%%%%%%%%%%%%%%%%%%%%%%%%%%%%%%%%%%%%%%%%%%%%%%%
\subsection{Linear dynamical systems}
In this section we investigate linear dynamical systems of the form
\begin{equation} \label{ode}
\begin{array}{rcl}  
E(p) \dot{x}(t,p) & = & 
A(p) x(t,p) + B(p) u(t) \\[1ex] 
y(t,p) & = & C(p) x(t,p) , \\
\end{array} 
\end{equation}
where the matrices $A,E \in \real^{n \times n}$, 
$B \in \real^{n \times n_{\rm in}}$ and 
$C \in \real^{n_{\rm out} \times n}$ depend on 
physical parameters $p \in \Pi \subseteq \real^{n_{\rm par}}$. 
An input $u : [0,\infty) \rightarrow \real^{n_{\rm in}}$ is given
and an output 
$y : [0,\infty) \times \Pi \rightarrow \real^{n_{\rm out}}$
is provided.
Without loss of generality, we restrict the investigations to the case of
single-input-single-output (SISO), i.e., $n_{\rm in} = n_{\rm out} = 1$.

We assume that the mass matrix~$E$ is always non-singular.
Consequently, a system of ODEs~(\ref{ode}) is given
with the state variables 
$x : [0,\infty) \times \Pi \rightarrow \real^n$. 
In our examinations,
initial value problems (IVPs) $x(0,p)=0$ are predetermined
for all~$p \in \Pi$.
Furthermore, we assume that the systems~(\ref{ode}) are asymptotically stable
for all $p \in \Pi$,
i.e., all eigenvalues $\Sigma(p) \subset \complex$ 
of the matrix pencil $\lambda E(p) - A(p)$ 
have a negative real part.  

The input-output mapping of the system~(\ref{ode}) can be
specified by a transfer function
$H : ( \complex \backslash \Sigma(p)) \rightarrow \complex$
in the frequency domain,
see~\cite[p.~65]{antoulas}. % Eq. (4.22)
The transfer function of the system~(\ref{ode}) becomes
\begin{equation} \label{transferfcn}
H(s,p) := C(p) \left( s E(p) - A(p) \right)^{-1} B(p) 
\qquad \mbox{for} \;\; s \in \complex \backslash \Sigma(p) ,
\end{equation}
which represents a rational function in the frequency variable~$s$.
The input-output mapping reads as $Y(s,p) = H(s,p) U(s)$,
where $U,Y$ denote the Laplace transforms of the input and the output,
respectively.

%%%%%%%%%%%%%%%%%%%%%%%%%%%%%%%%%%%%%%%%%%%%%%%%%%%%%%%%%%%%%%%%%%%%%%%%%%%%%
\subsection{Stochastic modeling and orthogonal basis}
In the dynamical system~(\ref{ode}), 
the parameters~$p \in \Pi$ are replaced by independent random variables
$p : \Omega \rightarrow \Pi$ 
on some probability space $(\Omega,\mathscr{A},\mu)$ 
with event space~$\Omega$, sigma-algebra~$\mathscr{A}$ and 
probability measure~$\mu$. 
We suppose the existence of a joint probability density function
$\rho : \Pi \rightarrow \real$.
Given a measurable function $f : \Pi \rightarrow \real$, the expected value 
reads as
$$ \mathbb{E} \left[ f \right] := 
\int_{\Omega} f (p(\omega)) \; \mbox{d}\mu(\omega) = 
\int_{\Pi} f (p) \rho(p) \; \mbox{d}p $$
provided that the integral is finite.
The associated Hilbert space
$$ \ltwo := \left\{ f : \Pi \rightarrow \real \; : \; 
f \; \mbox{measurable and} \;
\mathbb{E} \left[ f^2 \right] < \infty \right\} $$
features the inner product
\begin{equation} \label{innerproduct} 
\langle f , g \rangle := \mathbb{E} \left[ f g \right] =  
\int_{\Pi} f (p) g (p) \rho(p) \; \mbox{d}p 
\qquad \mbox{for} \;\; f,g \in \ltwo . 
\end{equation}
The norm of the Hilbert space reads as
\begin{equation} \label{L2-norm}
  \left\| f \right\|_{\ltwo} := \sqrt{ \langle f , f \rangle } .
\end{equation}
In the system~(\ref{ode}), we assume that 
$x_1(t,\cdot),\ldots,x_n(t,\cdot), y(t,\cdot) \in \ltwo$
point-wise for $t \in [0,\infty)$. 

We consider an orthonormal system $( \Phi_i )_{i \in \mathbbm{N}} \subset \ltwo$
consisting of polynomials.
The theory of the generalized polynomial chaos (gPC) applies
in this situation, see~\cite{xiu-karniadakis,xiu-book}. 
Hence the basis functions satisfy the orthogonality relations
\begin{equation} \label{orthogonal}
  \langle \Phi_i , \Phi_j \rangle = \left\{
  \begin{array}{ll}
    0 & \mbox{for} \;\; i \neq j , \\
    1 & \mbox{for} \;\; i = j . \\
  \end{array}
  \right.
\end{equation}
Such an orthonormal system exists under the above assumptions.
However, the system is not always complete, cf.~\cite{ernst}.
Since we investigate finite approximations, completeness is not
required in the following.

Finite approximations are obtained by
\begin{equation} \label{pce} 
x^{(m)}(t,p) = \sum_{i=1}^m v_i(t) \Phi_i(p) 
\qquad \mbox{and} \qquad 
y^{(m)}(t,p) = \sum_{i=1}^m w_i(t) \Phi_i(p) 
\end{equation}
including $m$ basis polynomials.
Due to the orthogonalities~(\ref{orthogonal}),
the transient coefficient functions 
$v_i : [0,\infty) \rightarrow \real^n$ and 
$w_i : [0,\infty) \rightarrow \real$
given by
\begin{equation} \label{coefficients}
v_{i,j}(t) = \langle x_j(t,\cdot) , \Phi_i(\cdot) \rangle 
\qquad \mbox{and} \qquad 
w_i(t) = \langle y(t,\cdot) , \Phi_i(\cdot) \rangle  
\end{equation}
represent the best approximation in the subspace spanned by
$\{\Phi_1,\ldots,\Phi_m\}$ with respect to the norm~(\ref{L2-norm}).
If the orthonormal system is complete, then the convergence
$$ \lim_{m \rightarrow \infty}
\left\| y(t,\cdot) - y^{(m)}(t,\cdot) \right\|_{\ltwo} = 0
\qquad \mbox{for each} \;\; t $$
is guaranteed for the output (QoI) and, likewise, for the state variables.

We assume that the orthonormal system $( \Phi_i )_{i \in \mathbbm{N}}$
is ordered in ascending total degrees, i.e.,
${\rm degree}(\Phi_i) \le {\rm degree}(\Phi_j)$ for $i \le j$.
Thus $\Phi_1 \equiv 1$ is the unique constant polynomial.
Let the approximation~(\ref{pce}) include all polynomials up to
a total degree~$d$.
The number of polynomials reads as, see~\cite[p.~65]{xiu-book},
% Eq. (5.24)
\begin{equation} \label{polynomial-number}
m(d) = \frac{(n_{\rm par}+d)!}{n_{\rm par}! d!} 
\end{equation}
for a fixed number~$n_{\rm par}$ of random parameters.
This number becomes large for high dimensions~$n_{\rm par}$ of the random space,
even though the total degree may be moderate, say $2 \le d \le 5$.

%%%%%%%%%%%%%%%%%%%%%%%%%%%%%%%%%%%%%%%%%%%%%%%%%%%%%%%%%%%%%%%%%%%%%%%%%%%%%
\subsection{Stochastic Galerkin method}

We investigate the input-output behavior of a
stochastic Galerkin system.
%%%%%%%%%%%%%%%%%%%%%%%%%%%%%%%%%%%%%%%%%%%%%%%%%%%%%%%%%%%%%%%%%%%%%%%%%%%%%
\subsubsection{Stochastic Galerkin system and transfer function}
\label{sec:galerkin} 
An approximation of the coefficient functions~(\ref{coefficients})
belonging to the solution of the random dynamical system~(\ref{ode})
can be obtained by a stochastic Galerkin approach.
This technique yields a larger coupled system
\begin{equation} \label{galerkin}
\begin{array}{rcl}  
\hat{E} \dot{\hat{v}}(t) & = & 
\hat{A} \hat{v}(t) + \hat{B} u(t) \\[1ex] 
\hat{w}(t) & = & \hat{C} \hat{v}(t) \\
\end{array} 
\end{equation}
with matrices $\hat{A},\hat{E} \in \real^{mn \times mn}$,
$\hat{B} \in \real^{mn}$ and $\hat{C} \in \real^{m \times mn}$.
Thus the system is single-input-multiple-output (SIMO).
The matrix~$\hat{E}$ is non-singular in most cases
and thus~(\ref{galerkin}) is a system of deterministic ODEs.
Although the stochastic Galerkin system may be unstable,
see~\cite{pulch-augustin}, this loss of stability happens rather seldom.
We assume that the system~(\ref{galerkin}) is asymptotically stable.
IVPs $\hat{v}(0)=0$ are imposed.
More details on the stochastic Galerkin approach for random linear dynamical
systems can be found in~\cite{pulch11,pulch14,pulch2018}.

The outputs $\hat{w} = (\hat{w}_1,\ldots,\hat{w}_m)^\top$
of the system~(\ref{galerkin}) represent
an approximation of the exact coefficients~(\ref{coefficients}).
The associated approximation of the QoI in the
random dynamical system~(\ref{ode}) reads as
\begin{equation} \label{qoi-appr}
  \hat{y}^{(m)}(t,p) = \sum_{i=1}^m \hat{w}_i(t) \Phi_i(p) .
\end{equation}
The Galerkin-projected system~(\ref{galerkin}) exhibits its own
input-output behavior described by a complex-valued transfer function
$\hat{H} : \complex \backslash \hat{\Sigma} \rightarrow \complex$,
$\hat{H} = (\hat{H}_1,\ldots,\hat{H}_m)^\top$
with a finite set of poles~$\hat{\Sigma}$,
see~\cite{pulch-maten-augustin,pulch-maten}.
It holds that
\begin{equation} \label{transfer-galerkin}
    \hat{H}(s) = \hat{C} \left( s \hat{E} - \hat{A} \right)^{-1} \hat{B}
    \qquad \mbox{for} \;\; s \in \complex \backslash \hat{\Sigma} .
\end{equation}
Again the transfer function represents a rational function.
The linear dynamical system~(\ref{galerkin}) is asymptotically stable,
if and only if the set of poles~$\hat{\Sigma}$ is located in the left half
of the complex plane. 
Moreover, the transfer function is always strictly proper for
systems of ODEs,
which is defined by the condition
\begin{equation} \label{transfer-to-zero}
  \lim_{s \rightarrow \infty} \hat{H}(s) = 0
\end{equation}
in each component.
The computational effort for an evaluation of~(\ref{transfer-galerkin})
is dominated by a single $LU$-decomposition of the matrix
$s\hat{E} - \hat{A}$.

In the frequency domain, Hardy norms characterize the magnitude
of a transfer function, see~\cite{antoulas}.
These norms allow for error estimates.
We will use the $\mathscr{H}_2$-norm component-wise, i.e., 
\begin{equation} \label{hardy-norm}
  \left\| \hat{H}_i \right\|_{\mathscr{H}_2} =
  \sqrt{ \displaystyle \frac{1}{2\pi} \int_{-\infty}^{+\infty} 
    \left| \hat{H}_i ({\rm i} \, \omega) \right|^2 \; {\rm d}\omega }
  \; =
  \sqrt{ \displaystyle \frac{1}{\pi} \int_0^{+\infty} 
    \left| \hat{H}_i ({\rm i} \, \omega) \right|^2 \; {\rm d}\omega } 
\end{equation}
for $i=1,\ldots,m$ with ${\rm i} = \sqrt{-1}$.
Alternatively, the $\mathscr{H}_{\infty}$-norm can be
considered component-wise, i.e.,
\begin{equation} \label{hinf-norm}
  \left\| \hat{H}_i \right\|_{\mathscr{H}_{\infty}} =
  \sup_{\omega \in \real} 
    \left| \hat{H}_i ({\rm i} \, \omega) \right|
  =
   \sup_{\omega \ge 0} 
    \left| \hat{H}_i ({\rm i} \, \omega) \right|
\end{equation}
for $i=1,\ldots,m$.
In~\cite{pulch2018}, both $\mathscr{H}_2$- and $\mathscr{H}_{\infty}$-norm
were investigated in this context.
We apply only the $\mathscr{H}_2$-norm in the following,
because the results are qualitatively the same in both cases.

%%%%%%%%%%%%%%%%%%%%%%%%%%%%%%%%%%%%%%%%%%%%%%%%%%%%%%%%%%%%%%%%%%%%%%%%%%%%%
\subsubsection{Computation of Hardy norms}
\label{sec:quadrature}
There are several possibilities to compute the $\mathscr{H}_2$-norm
of a linear dynamical system, cf.~\cite{antoulas}.
We use an own technique based on quadrature, where the computational effort
is nearly independent of the number~$m$ of basis polynomials.
An approximation of the $\mathscr{H}_2$-norms~(\ref{hardy-norm}) is obtained
by a quadrature with (single) rectangular rule in $[0,\omega_{\min}]$ and
composite trapezoidal rule in $[\omega_{\min},\omega_{\max}]$ given
$0<\omega_{\min}<\omega_{\max}$.
We apply a grid of the form
$\omega_{\min} = \omega_1 < \omega_2 < \cdots < \omega_{\nu-1} < \omega_{\nu}=\omega_{\max}$.
% Let $\Delta \omega_j = \omega_j - \omega_{j-1}$ for $j=2,\ldots,\nu$.
The approximation of~(\ref{hardy-norm}) reads as
$$  \int_0^{\infty} \left| \hat{H}_i({\rm i}\omega) \right|^2 \; {\rm d}\omega
\approx
% \left| \hat{H}_i ({\rm i} \, \omega) \right|^2 \; {\rm d}\omega \approx
\omega_1 \left| \hat{H}_i({\rm i}\omega_1) \right|^2 +
\sum_{j=2}^{\nu} \frac{\omega_j - \omega_{j-1}}{2}
\left( \left| \hat{H}_i({\rm i}\omega_{j-1}) \right|^2
+ \left| \hat{H}_i({\rm i}\omega_j) \right|^2 \right) $$
for $i=1,\ldots,m$.
Since the transfer function~(\ref{transfer-galerkin}) is continuous
at~$\omega=0$, the quadrature error within $[0,\omega_{\min}]$
is negligible for $\omega_{\rm \min}$ sufficiently close to zero.
The property~(\ref{transfer-to-zero}) guarantees that the truncation error,
which arises by discarding the interval $(\omega_{\max},\infty)$,
is negligible for sufficiently large~$\omega_{\max}$.
In each node of the quadrature, the evaluation of~(\ref{transfer-galerkin})
requires mainly an $LU$-decomposition of the matrix
${\rm i}\omega_j \hat{E} - \hat{A}$
independent of the output matrix~$\hat{C}$ provided that~$m$
is not extremely large.

%%%%%%%%%%%%%%%%%%%%%%%%%%%%%%%%%%%%%%%%%%%%%%%%%%%%%%%%%%%%%%%%%%%%%%%%%%%%%
\subsubsection{Applicability to differential-algebraic equations}
If the mass matrix~$E(p)$ is singular, then the linear dynamical
system~(\ref{ode}) consists of differential-algebraic equations (DAEs).
In most cases, the mass matrix~$\hat{E}$ of the stochastic Galerkin
system~(\ref{galerkin}) also becomes singular.
We assume an asymptotically stable system~(\ref{galerkin}) again,
where the set of poles~$\hat{\Sigma}$ is situated in the left half
of the complex plane.
It follows that the associated matrix pencil is regular.
%A regular matrix pencil is assumed, i.e.,
%$\det (s \hat{E} - \hat{A}) \neq 0$ for almost
%all $s \in \complex$.
We outline the potential to obtain the frequency domain information
from Section~\ref{sec:galerkin}, which is required for the
basis selection in Section~\ref{sec:fd-analysis}.

A linear system of DAEs is characterized by its (nilpotency) index
$\kappa \in \nat$ ($\kappa \ge 1$), see~\cite[p.~454]{hairer2}.
Let $\kappa$ be the index of the stochastic Galerkin system~(\ref{galerkin}).
Consequently, the transfer function of~(\ref{galerkin}) reads as
\begin{equation} \label{transfer-dae}
  \hat{H}(s) = \hat{H}_{\rm SP}(s) + \hat{P}(s)
\end{equation}
with a strictly proper vector-valued rational function~$\hat{H}_{\rm SP}$
satisfying~(\ref{transfer-to-zero}) and
a vector-valued polynomial~$\hat{P}$ of degree at most~$\kappa-1$.
The formula~(\ref{transfer-dae}) is valid for DAEs in general,
see~\cite{benner-stykel}.
If the polynomial part~$\hat{P}$ vanishes, then the
$\mathscr{H}_2$-norms~(\ref{hardy-norm}) do exist for all components.
Hence the strategy performs as in the case of ODEs.
If the polynomial part does not vanish, then some or all
$\mathscr{H}_2$-norms do not exist.

If the index of the stochastic Galerkin system~(\ref{galerkin}) is one,
then the polynomial part~$\hat{P}$ becomes constant (degree zero).
In this case, the existence of the $\mathscr{H}_{\infty}$-norms~(\ref{hinf-norm})
is guaranteed for all components.
The frequency domain analysis can be done using the norms~(\ref{hinf-norm}).

%%%%%%%%%%%%%%%%%%%%%%%%%%%%%%%%%%%%%%%%%%%%%%%%%%%%%%%%%%%%%%%%%%%%%%%%%%%%%
\subsection{Basis selection and low-dimensional representations}
Let an initial approximation~(\ref{pce}) be given with a large
number~(\ref{polynomial-number}) of basis polynomials up to total degree~$d$.
Let $\mathcal{I}^d := \{ 1,\ldots,m(d) \}$ be the respective index set.
In general, nearly all associated coefficient functions are non-zero.
However, the coefficients typically exhibit different orders of magnitudes
and a fast decay for increasing degree.

Our aim is to identify a low-dimensional representation
\begin{equation} \label{low-dim-appr}
  \tilde{y}^{(\mathcal{I})} (t,p) = \sum_{i \in \mathcal{I}}
  \tilde{w}_i(t) \Phi_i(p)
\end{equation}
with an index set $\mathcal{I} \subset \mathcal{I}^d$
satisfying $q := | \mathcal{I} | \ll m(d)$.
Thus the selected basis functions from~(\ref{low-dim-appr})
represent a subset of the basis functions in~(\ref{pce}).
Yet the difference $y - \tilde{y}^{(\mathcal{I})}$ with $y$ from~(\ref{ode})
or (at least) $\hat{y}^{(m(d))} - \tilde{y}^{(\mathcal{I})}$
with $\hat{y}$ from~(\ref{qoi-appr})
should be sufficiently
small in the $\ltwo$-norm point-wise in time.
The approximation~(\ref{low-dim-appr}) is also called a
$q$-sparse representation, cf.~\cite{blatman} for a motivation.
Hence the problem consists in the identification of both
the index set~$\mathcal{I}$ and associated approximations~$\tilde{w}_i$
of the coefficient functions for a desired number~$q$.
If a transient solution of the stochastic Galerkin system~(\ref{galerkin})
is available, then we can simply set $\tilde{w}_i := \hat{w}_i$
for $i \in \mathcal{I}$.

%\begin{equation} \label{index-large}
%\mathcal{I}^d := \left\{ i \; : \; 
%\Phi_i(p) = \phi_{j_1}^{(1)}(p_1) \cdots \phi_{j_q}^{(q)}(p_q) 
%\;\; \mbox{with} \;\; j_1 + \cdots + j_q \le d 
%\right\} . 
%\end{equation}
% systems $(\phi_j^{(\ell)})_{j \in \mathbbm{N}_0}$ 

%%%%%%%%%%%%%%%%%%%%%%%%%%%%%%%%%%%%%%%%%%%%%%%%%%%%%%%%%%%%%%%%%%%%%%%%%%%%%
\subsection{Analysis in frequency domain}
\label{sec:fd-analysis}
In~\cite{pulch2018}, the error $\hat{y}^{(m(d))} - \tilde{y}^{(\mathcal{I})}$
between the approximation~(\ref{qoi-appr}) from the stochastic Galerkin
method and the low-dimensional approximation~(\ref{low-dim-appr}) was
analyzed using the frequency domain.
The $\mathscr{H}_2$-norms~(\ref{hardy-norm}) of the transfer function
belonging to the Galerkin-projected system~(\ref{galerkin}) are considered.
Let $\mathcal{I} \subset \mathcal{I}^d$ be any index set. 
Theorem~1 in~\cite{pulch2018} demonstrates the error estimate
\begin{equation} \label{errorbound} 
\sup_{t \ge 0} \left\| \hat{y}^{(m(d))} (t,\cdot) - 
\tilde{y}^{(\mathcal{I})}(t,\cdot) \right\|_{\ltwo} \le 
\sqrt{ \displaystyle \sum_{i \in \mathcal{I}^d \backslash \mathcal{I}} 
\left\| \hat{H}_i \right\|_{\mathscr{H}_2}^2 } \;\;
\| u \|_{\ltwotime} 
\end{equation}
with $\tilde{w}_i = \hat{w}_i$ for all $i \in \mathcal{I}$.
The usual $\ltwotime$-norm is used for the time-dependent input.
An advantage is that the error is bounded uniformly for all times.
However, this error bound is not sharp.

For a desired cardinality~$q = |\mathcal{I}|$,
we determine a minimum upper bound in the right-hand side
of~(\ref{errorbound}).
An optimal index set satisfies, see~\cite[p.~7]{pulch2018},
\begin{equation} \label{optimal-index-set}
  \mathcal{I}_q = \underset{\mathcal{I} \subseteq \mathcal{I}^d}{\arg\min}
  \left\{ \displaystyle \sum_{i \in \mathcal{I}^d \backslash \mathcal{I}} 
  \left\| \hat{H}_i  \right\|_{\mathscr{H}_2}^2 :
  \left| \mathcal{I} \right| = q \right\} .
\end{equation}
In practice, the norms will be pairwise different.
Hence the unique index set~$\mathcal{I}_q$ includes the
components with the $q$~largest $\mathscr{H}_2$-norms.
Furthermore, the minimization of the upper bound is independent of the
choice of the input~$u$.
Approximations of the $\mathscr{H}_2$-norms,
which occur in~(\ref{optimal-index-set}), can be computed numerically
as described in Section~\ref{sec:quadrature}. 

This frequency domain analysis can be extended
to the case of multiple outputs ($n_{\rm out} > 1$) in~(\ref{ode})
straightforward.
Either a large output matrix $\hat{C} \in \real^{mn_{\rm out} \times mn}$
in~(\ref{galerkin}) or $n_{\rm out}$ separate systems~(\ref{galerkin})
with output matrices
$\hat{C}_j \in \real^{m \times mn}$ for $j=1,\ldots,n_{\rm out}$ can be arranged.
We discuss the latter approach.
In the transfer function~(\ref{transfer-galerkin}),
the part
$$ \hat{F}(s) := (s \hat{E} - \hat{A})^{-1} \hat{B} $$
is identical for any output.
An evaluation of $\hat{F}({\rm i}\omega)$ can be reused for all systems
in a quadrature to approximate the integrals~(\ref{hardy-norm}).
Thus the computational effort for a frequency domain analysis is nearly
independent of the number of outputs (provided that $n_{\rm out} \ll mn$).
A separate low-dimensional basis can be computed for each output.

%%%%%%%%%%%%%%%%%%%%%%%%%%%%%%%%%%%%%%%%%%%%%%%%%%%%%%%%%%%%%%%%%%%%%%%%%%%%%

%%%%%%%%%%%%%%%%%%%%%%%%%%%%%%%%%%%%%%%%%%%%%%%%%%%%%%%%%%%%%%%%%%%%%%%%%%%%%
%%%                       Numerical Methods                               %%%
%%%%%%%%%%%%%%%%%%%%%%%%%%%%%%%%%%%%%%%%%%%%%%%%%%%%%%%%%%%%%%%%%%%%%%%%%%%%%

\section{Numerical Methods}
\label{sec:methods}
In this section, we discuss numerical techniques for a computation of
a low-dimensional representation and its coefficient functions
in the time domain. Spe\-ci\-fi\-cally, we consider two approaches: (i) sparsity inducing regression and (ii) least squares regression using the sparse index set consisting of the $q$~largest Hardy norms identified by the frequency domain analysis.

%%%%%%%%%%%%%%%%%%%%%%%%%%%%%%%%%%%%%%%%%%%%%%%%%%%%%%%%%%%%%%%%%%%%%%%%%%%%%
\subsection{Sampling techniques}
Our numerical methods apply the information from
$k$~samples of the output (QoI) of the
linear dynamical system~(\ref{ode}) for realizations
$p_1,\ldots,p_k$ of the random variables.
We collect the samples in a vector-valued time-dependent function
\begin{equation} \label{samples}
  \bar{y}(t) := ( y(t,p_1) , \ldots , y(t,p_k) )^\top .
\end{equation}
Let a basis $\{ \Phi_1 , \ldots , \Phi_{m(d)} \}$ be given with all
orthonormal polynomials up to a total degree~$d$.
We arrange the Vandermonde matrix
\begin{equation} \label{vandermonde}
  V \in \real^{k \times m} , \quad
  V = (v_{ij}) , \quad
  v_{ij} := \Phi_j(p_i) 
\end{equation}
for $i=1,\ldots,k$ and $j=1,\ldots,m$.
This matrix is dense.
% This matrix exhibits full column rank in the case of $k \ge m$.
Now let $v_1,\ldots,v_m \in \real^k$ be the columns of the
matrix~(\ref{vandermonde}) and
$\mathcal{I}_q \subseteq \mathcal{I}^d$ with
$\mathcal{I}_q = \{ j_1 , \ldots , j_q \}$ 
be an arbitrary subset of indices.
We use the notation
\begin{equation} \label{vandermonde-subset}
  V_{\mathcal{I}_{q}} := \left( v_{j_1} , \ldots , v_{j_q} \right)
  \in \real^{k \times q}
\end{equation}
for the matrix, which consists of the associated subset of columns
from the Vandermonde matrix~(\ref{vandermonde}).
The cardinality~\eqref{polynomial-number} of the total-degree basis $m(d)$ grows quickly with the number of random parameters $n_{\rm par}$ and degree~$d$. Thus, in the following, we assume $k < m(d)$, because for computationally expensive linear dynamical systems we cannot afford to run the model exhaustively.

%%%%%%%%%%%%%%%%%%%%%%%%%%%%%%%%%%%%%%%%%%%%%%%%%%%%%%%%%%%%%%%%%%%%%%%%%%%%%
\subsection{Orthogonal matching pursuit}
\label{sec:l0min}
We consider the determination of a sparse representation in the time domain now.

\subsubsection{Algorithm}
Let a fixed time point~$t > 0$ be given.
We use OMP~\cite{Pati1993} to determine a $q$-sparse representation
at this time point by greedily minimizing the $\ell_0$-norm
$$ \| w \|_0 := | \{ i \; : \; w_i \neq 0 \} |
\qquad \mbox{for} \;\; w \in \real^m , $$
which counts the number of non-zero entries in a vector.
Specifically, the OMP approach finds an approximation to the problem 
\begin{equation}
\label{eq:omp}
\tilde{w}(t) = \arg\min\; \| \tilde{w}(t) \|_{0}
\quad \text{such that}\quad
\|V\tilde{w}(t) - \bar{y}(t)\|_{2} \le \varepsilon
\end{equation}
with the Euclidean norm $\| \cdot \|_2$ and $\varepsilon \ge 0$ 
by iteratively building up an approximation of the solution vector $\tilde{w}(t)$. At each iteration step, a least squares problem is solved using a subset of active columns of $V$. At the $q$th iteration step, OMP updates an active index set $\mathcal{I}_q$, in a greedy fashion, such that the inactive column index $j\notin \mathcal{I}_q$ with the highest correlation (inner product) with the current residual is added to the active set.
This update reads as
\begin{align}\label{eq:omp-update}\mathcal{I}_q=\mathcal{I}_{q-1}\cup\{j\}
\quad \mbox{with} \quad
j = \operatornamewithlimits{\arg\max}_{i\in\mathcal{I}\setminus \mathcal{I}_{q-1}}\frac{r_{q-1}(t)^\top v_i}{\|v_i \|_2}.\end{align}
OMP then solves the least squares problem
\begin{equation}\label{omp-lstsq-step}
  \operatornamewithlimits{\arg\min}_{\tilde{w}^{(q)}(t)\in \real^q}\|r_q(t)\|_2 , 
\end{equation}
where $r_{q}(t)=V_{\mathcal{I}_q}\tilde{w}^{(q)}(t)-\bar{y}(t)$.

Setting $\varepsilon=0$ in~(\ref{eq:omp}),
the algorithm yields a finite sequence of
$q$-sparse approximations~(\ref{low-dim-appr}) for $q=1,2,\ldots,q_{\max}$, where $q_{\max}\le\min(k,m)$, thus removing the need to estimate $\varepsilon$.
In the following, we will assume that $q_{\max} \le k < m$, which always implies over-determined linear systems.
For a given number~$q$ the solution obtained using OMP is a greedy attempt to find a solution with $q$~non-zero coefficients
such that the residual
\begin{equation} \label{residual}
  R(t) := \| V \tilde{w}(t) - \bar{y}(t) \|_2
\end{equation}
is minimal. We will employ this approach in the time domain.

Unlike the least squares approach in Section~\ref{sec:lsp},
OMP must be applied at each time point separately. Here we use an OMP implementation based upon rank-one updates of a $QR$-factorization of the matrix $V_{\mathcal{I}_q}$ to form  $V_{\mathcal{I}_{q+1}}$, see~\cite[p.~334]{golub-loan}.
The total complexity of the $q$th iteration is $O(km+ks)$.
As shown in~\cite{strum2012},
the total complexity of our OMP for a single time point is
$$ q_{\max}km + m \sum_{q=1}^{q_{\max}} q =
m q_{\max} \left( k + \textstyle{\frac{1}{2}} (q_{\max}+1) \right). $$
OMP does not require the construction of the full Vandermonde matrix. If memory is limited, one can simply compute each basis vector $v_i$ for $i\in\mathcal{I}\setminus \mathcal{I}_{q-1}$ in each iteration $q=1,\ldots,q_{\max}$. The construction of the basis vector $v_i$ is necessary to compute the inner products $r_{q-1}(t)^\top v_i$ and $\lVert v_i\rVert_2$ in \eqref{eq:omp-update}. Although the full Vandermonde matrix does not have to be stored, the $QR$-factorization of the matrix $V_{\mathcal{I}_q}$ must be stored. Only if $q=m$, then $V_{\mathcal{I}_q}$ will be the full Vandermonde matrix. However, we suppose $q \ll m$.

  Note that avoiding the storage of the full Vandermonde matrix requires repeated computation of the basis vectors~$v_i$, which affects the complexity analysis. For the problem sizes considered in this paper memory was not an issue, so the full Vandermonde matrix was stored to avoid the redundant computation of the basis vectors. 

\subsubsection{Theoretical guarantees}
OMP is directly related to $\ell_1$-minimization~\cite{candes2006} which finds the dominant gPC coefficients by solving
\begin{equation}
\label{eq:bpdn}
\tilde{w}(t) = \arg\min\; \| \tilde{w}(t) \|_{1}
\quad \text{such that}\quad
\|V\tilde{w}(t) - \bar{y}(t)\|_{2} \le \varepsilon
\end{equation}
with the norm $\| w \|_1 = | w_1 | + \cdots + | w_m |$ for $w \in \real^m$.
This $\ell_1$-minimization problem is often referred to as
Basis Pursuit Denoising.
The problem obtained by setting $\varepsilon=0$, to enforce interpolation,
is termed Basis Pursuit (BP).  

Let $\{\Phi_1,\ldots,\Phi_m\}$ be an orthonormal polynomial basis
  with respect to the probability measure of the random variables $p$,
  i.e., the relations~(\ref{orthogonal}) are satisfied.
  Furthermore, in~(\ref{eq:bpdn}),
  let $V\in\real^{k\times m}$ be the matrix obtained by random sampling the
  basis under the probability measure. If the number of samples $k$ satisfies
$$
  k \ge C L_m q \log^3 q \log m
  \qquad \mbox{with} \qquad
  L_m = \max_{i=1,\ldots,m} \lVert\Phi_i\rVert^2_\infty,
$$
then every $q$-sparse vector can be recovered by
\eqref{eq:bpdn} (for $\varepsilon=0$) with probability at least
$1-m^{-\gamma\log^3 q}$
including constants $C,\gamma > 0$,
see~\cite{Rauhut_TFNMSR_2010}.
The maximum norm $\lVert\cdot\rVert_\infty$ is taken over
the support of the probability measure.

OMP requires stronger theoretical conditions than BP, see~\cite{cai2011}. 
For example, for a sufficiently small ratio of the number of samples to sparsity, BP can guarantee recovery of all $q$-sparse polynomials with high probability,
whereas, for a fixed set of samples, OMP can guarantee recovery of at least one sparse polynomial but not all, see~\cite{Kunis_R_FCM_2008}.
However, despite this theoretical difference, in practice OMP can still obtain
comparable accuracy to non-approximate algorithms such as BP~\cite{Kunis_R_FCM_2008}. Moreover, OMP has a much faster execution speed than most algorithms, which in conjunction with its iterative nature makes OMP more amenable to assessing the effect of sparsity on the accuracy of our approximations for dynamical systems.

We remark that although we focus on linear problems, sparse approximation strategies can also be applied to nonlinear dynamical systems, for example~\cite{hampton-doostan}.

%%%%%%%%%%%%%%%%%%%%%%%%%%%%%%%%%%%%%%%%%%%%%%%%%%%%%%%%%%%%%%%%%%%%%%%%%%%%%
\subsection{Least squares problem}
\label{sec:lsp}
Orthogonal matching pursuit greedily builds up a sparse representation of the output (QoI) of a linear dynamical system~(\ref{ode}) and uses least squares to solve an over-determined linear system. Alternatively, we determine a $q$-sparse representation applying a basis specified by the frequency domain analysis
in Section~\ref{sec:fd-analysis}, i.e.,
for the index set~(\ref{optimal-index-set}).

Let $V_{\mathcal{I}_q}$ be the matrix~(\ref{vandermonde-subset}) for
the index set~$\mathcal{I}_q$.
Likewise, let $\tilde{w}^{(q)}(t) \in \real^q$ be approximations
of the coefficient functions~(\ref{coefficients}) associated to
basis functions determined by the index set~$\mathcal{I}_q$.
Assuming $q < k$, we achieve a least squares problem 
\begin{equation} \label{least-squares}
  \min_{\tilde{w}^{(q)}(t) \in \real^q}
  \left\| V_{\mathcal{I}_q} \tilde{w}^{(q)}(t) - \bar{y}(t) \right\|_2
\end{equation}
with the Euclidean norm~$\| \cdot \|_2$
point-wise in time.
Since the matrix $V_{\mathcal{I}_q}$ exhibits full rank,
a unique solution $\breve{w}^{(q)}$ exists.
The solution $\breve{w}^{(q)}$ directly yields a
low-dimen\-sional approximation~(\ref{low-dim-appr}).

Moreover, the matrix~$V_{\mathcal{I}_q}$ is time-invariant and thus a
$QR$-factorization, see~\cite[p.~246]{golub-loan}, % Ch. 5.2
of this matrix can be reused to solve the least squares problems
in all time points.
% QR-decomp: Ch. 4.7 in Stoer/Bulirsch
This property makes the technique also advantageous in the
case of $q \approx k$ ($q<k$), if $k$ is not too large.
Using Householder transformations, the number of operations
becomes $2kq^2 - \frac{2}{3}q^3$ in the $QR$-factorization.
The numerical solution of the problem~(\ref{least-squares}) with a given
$QR$-decomposition requires just about
$2kq + \frac{1}{2}q^2$ operations.
% 2kq for Q^T y and q^2/2 for backward substitution
Proceeding for $q=1,\ldots,q_{\max}$,
rank-one updates can be used in the $QR$-decompositions,
see~\cite[p.~334]{golub-loan}, % Ch. 12.5
for keeping the computational complexity small.

\subsection{Related work}
  The goal of this paper is to explore methods for producing sparse representations of the output of linear dynamical systems. Specifically, we consider time-frequency analysis and compressive sensing to produce sparse representations. Other approaches do exist for reducing the complexity of a polynomial approximation of a function. For example, low-rank approaches~\cite{Doostan_2013,Oseledets_2013} can be used to find low-rank approximations of the coefficients of a tensor-product polynomial basis. For tensor-train representations~\cite{Gorodetsky_2018b}, the number of unknowns grows only linearly with dimension and quadratically with rank. The coefficients of a tensor-product basis are often dense. However, there have been some attempts to combine low-rank approximation with sparsity inducing methods~\cite{Mathelin_2014,Zhang_WD_IEEETCPMT_2017}.

  The cardinality~\eqref{polynomial-number} of the total-degree polynomial basis grows exponentially with dimension $d$ and, consequently, can become intractable for large number of variables. Often only a subset of variables and interactions between variables contribute significantly to variations in a function. The (adaptive) ANOVA decomposition can be used to identify and exploit such effective dimensionality~\cite{Prasad_R_IEEETMTT_2017,Yang_CLK_JCP_2012,Zhang_YOKD_IEEECDICS_2015}. A method that exploits both sparsity and adaptive ANOVA was developed in~\cite{jakeman2015}.

  The ability to recover sparse, and even low-rank approximations, is dependent on the method used to select samples at which to evaluate the function being approximated. For sparse approximation using Legendre basis randomly sampling from the uniform probability measure, as we do in Section~\ref{sec:example}, is adequate. However, for other orthogonal polynomials sampling from the probability measure can dramatically decrease the ability to recover sparse approximations. Changes of measure have been used successfully to improve sparse recovery~\cite{Hampton_D_JCP_2015,Jakeman_NZ_SISC_2017,Rauhut_W_JAT_2012}. Algorithms for generating samples with a well-conditioned Vandermonde matrix via subsampling of tensor-product quadrature have also been used for sparse approximation~\cite{Tang_2014} as well as regression, interpolation and low-rank approximation~\cite{Ahadi_R_IEEETCDICS_2016,Li_Z_WRR_2007,Manfredi_GZC_IEEEMWCL_2015,Zhang_MED_IEEETCDICS_2013}.

%%%%%%%%%%%%%%%%%%%%%%%%%%%%%%%%%%%%%%%%%%%%%%%%%%%%%%%%%%%%%%%%%%%%%%%%%%%%%
\subsection{Error measures in time domain}
We assume that the solution of the stochastic Galerkin
system~(\ref{galerkin}) yields a sufficiently accurate
approximation~(\ref{qoi-appr}) of the QoI
by the selection of a sufficiently large total degree~$d$.
Thus the truncation error as well as the error of the Galerkin approach are
negligible in comparison to a sparsification error.
A low-dimen\-sional approximation~(\ref{low-dim-appr}) is defined by
the index set~$\mathcal{I}$ and its coefficient functions
$\tilde{w}_i$ for $i \in \mathcal{I}$.
We embed these coefficients in a vector $\tilde{w} \in \real^{m}$
by specifying $\tilde{w}_i = 0$ for $i \notin \mathcal{I}$.

In the time domain, the $\ltwo$-norm~(\ref{L2-norm}) of the error
becomes
$$ \left\| \hat{y}^{(m(d))} (t,\cdot) - \tilde{y}^{(\mathcal{I})} (t,\cdot)
\right\|_{\ltwo} = \left\| \hat{w}(t) - \tilde{w}(t) \right\|_2
\qquad \mbox{for} \;\; t \ge 0 $$
with the Euclidean norm $\| \cdot \|_2$ due to Parseval's theorem.
We examine the relative error in $\ltwo$, i.e.,
\begin{equation} \label{l2error}
  E_{\mathscr{L}^2}(t) :=
  \frac{\| \hat{w}(t) - \tilde{w}(t) \|_2}{\| \hat{w}(t) \|_2}
\end{equation}
point-wise for~$t \ge 0$.

The identification of low-dimensional representations
uses $k$ samples~(\ref{samples}) of the QoI.
% The residual of a low-dimensional approximation~(\ref{low-dim-appr})
% becomes 
% \begin{equation} \label{residual}
%   R(t) := \| V \tilde{w}(t) - \bar{y}(t) \|_2
% \end{equation}
% for $t \ge 0$ with the Vandermonde matrix~(\ref{vandermonde}).
Both the  minimization problem from Section~\ref{sec:l0min} and the
least squares problem from Section~\ref{sec:lsp} feature the
property that the residual~(\ref{residual}) decreases monotone
for increasing dimension~$q$ of the subspace at fixed time.
The residual~(\ref{residual}) also yields a rough estimate of
the $\ltwo$-norm of the difference between the exact
QoI and the low-dimensional representation, i.e.,
\begin{equation} \label{L2-estimate}
  \left\| y (t,\cdot) - \tilde{y}^{(\mathcal{I})} (t,\cdot) \right\|_{\ltwo}
  \approx \textstyle
  \frac{1}{\sqrt{k}} \| V \tilde{w}(t) - \bar{y}(t) \|_2
\end{equation}
for each~$t$.
The estimate~(\ref{L2-estimate}) converges for $k,m \rightarrow \infty$
in probability (provided that the orthonormal basis is complete).

The error measures~(\ref{l2error}) and~(\ref{residual}) are given
point-wise in time.
Global measures can be obtained by taking (integral) mean values
on a finite time interval.
In practice, the problems are solved for a finite set of time points
and (arithmetic) means can be applied.

%%%%%%%%%%%%%%%%%%%%%%%%%%%%%%%%%%%%%%%%%%%%%%%%%%%%%%%%%%%%%%%%%%%%%%%%%%%%%
%%%                         Test Example                                  %%%
%%%%%%%%%%%%%%%%%%%%%%%%%%%%%%%%%%%%%%%%%%%%%%%%%%%%%%%%%%%%%%%%%%%%%%%%%%%%%

\section{Test Example}
\label{sec:example}
\begin{figure}[t]
\begin{center}
\includegraphics[width=5cm]{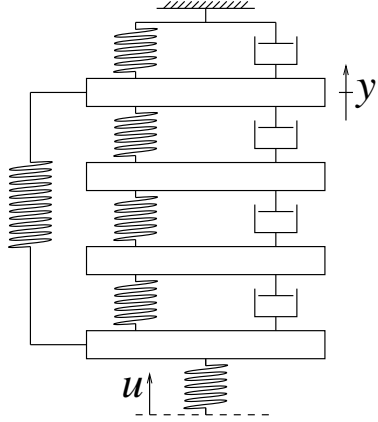}
\end{center}
\caption{Mass-spring-damper configuration.}
\label{fig:mass-spring-damper}
\end{figure}

We use a mass-spring-damper system from~\cite{lohmann-eid}
depicted in Figure~\ref{fig:mass-spring-damper}.
This configuration consists of 4~masses, 6~springs and 4~dampers,
i.e., $n_{\rm par}=14$ parameters.
The input is the excitation at the bottom mass, whereas the position
of the top mass represents the output.
The mathematical modeling yields a system~(\ref{ode}) of ODEs with 
dimension $n=8$.

%%% Figure: Statistics %%%%%%%%%%%%%%%%%%%%%%%%%%%%%%%%%%%%%%%%%%%%%%%%%%%%%%
%%%%%%%%%%%%%%%%%%%%%%%%%%%%%%%%%%%%%%%%%%%%%%%%%%%%%%%%%%%%%%%%%%%%%%%%%%%%%

Transient simulations are performed in a time interval $[0,T]$.
We supply a harmonic oscillation as input signal 
\begin{equation} \label{input-periodic}
  u(t) = \left\{ \begin{array}{cl}
    \sin ( \omega_0 t ) & \mbox{for} \;\; 0 \le t \le T , \\
    0 & \mbox{otherwise} , \\
  \end{array} \right.
\end{equation}
with frequency $\omega_0$.
The compact support $[0,T]$ guarantees $u \in \ltwotime$.
We choose $T=500$ and $\omega_0 = 0.1$ in this test example.
It holds that $\| u \|_{\ltwotime} \approx 16$.
% sqrt( T/2 ) = 15.8

In the following, all IVPs of ODEs are solved numerically
by a Runge-Kutta method of order~4(5)
with a step size selection based on the relative tolerance
$\varepsilon_{\rm rel} = 10^{-6}$ and the absolute tolerance
$\varepsilon_{\rm abs} = 10^{-8}$.
These high accuracy requirements imply that the errors of the
time integrations are negligible in comparison to the other
error sources.

In the stochastic modeling, we choose the means of parameters
as the values given in~\cite{lohmann-eid}.
Independent uniformly distributed random variables are introduced
with 5\% variation around the mean values.
Thus the gPC expansions~(\ref{pce}) include the multivariate
Legendre polynomials.
We incorporate all polynomials up to total degree $d=3$,
which results in $m=680$ basis functions:
1 of degree zero, 14 of degree one, 105 of degree two and
560 of degree three.

%%% Figure: Statistics %%%%%%%%%%%%%%%%%%%%%%%%%%%%%%%%%%%%%%%%%%%%%%%%%%%%%%
\begin{figure}
\begin{center}
\includegraphics[width=7cm]{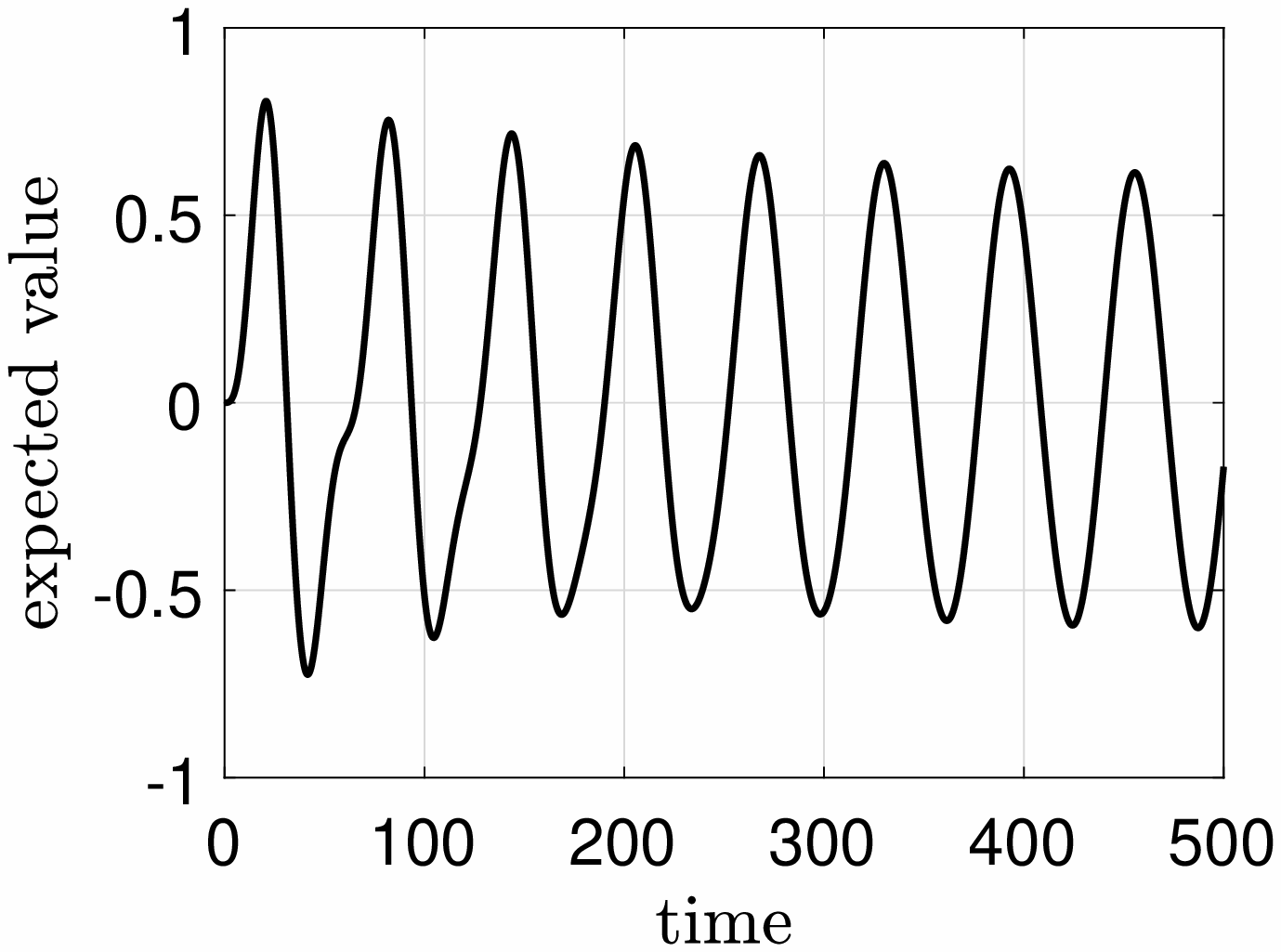}
\hspace{10mm}
\includegraphics[width=7cm]{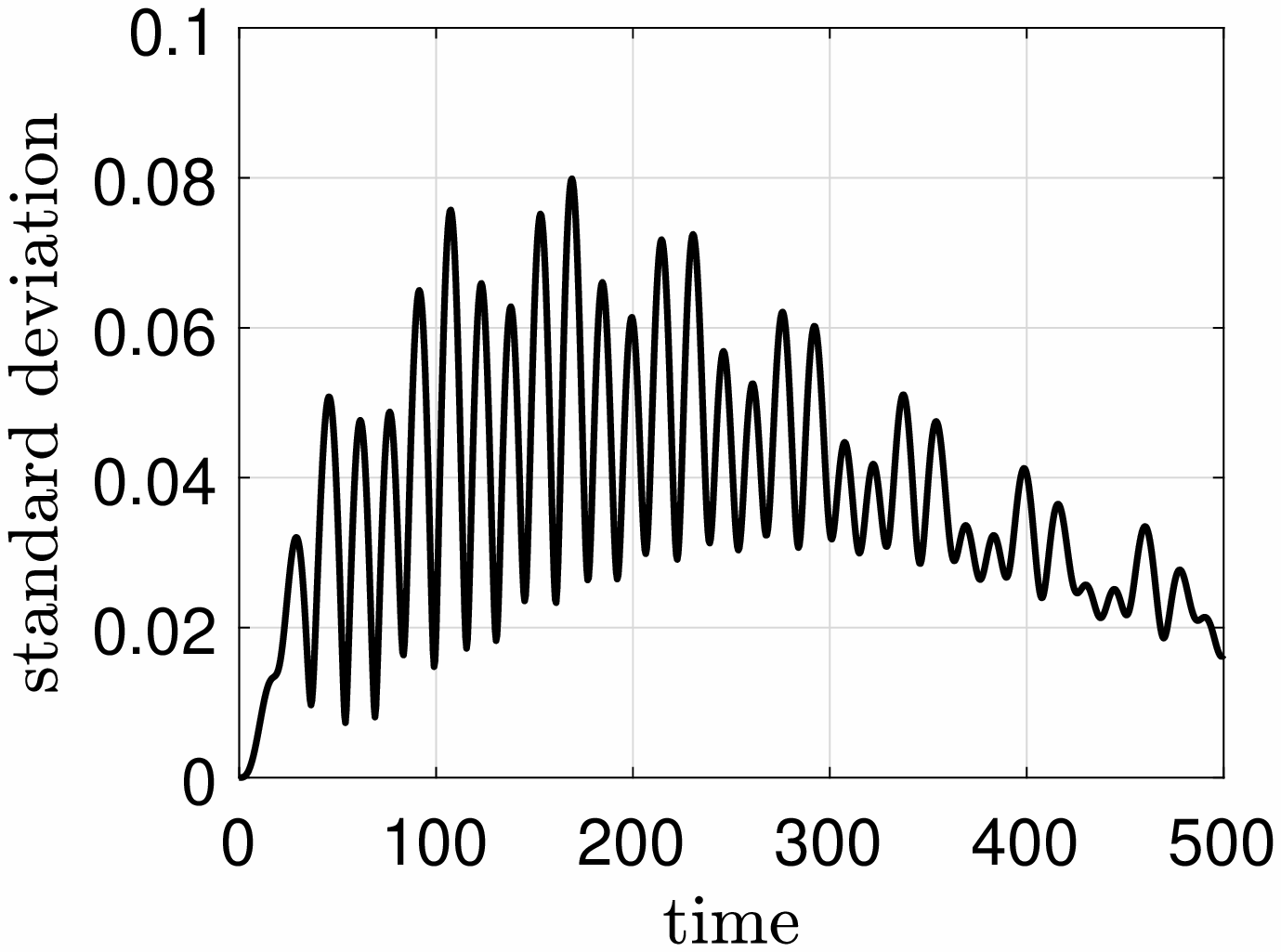}
\end{center}
\caption{Expected value (left) and standard deviation (right) of
  random output for periodic input signal in mass-spring-damper system.}
\label{fig:statistics}
\end{figure}
%%%%%%%%%%%%%%%%%%%%%%%%%%%%%%%%%%%%%%%%%%%%%%%%%%%%%%%%%%%%%%%%%%%%%%%%%%%%%

Now we arrange the stochastic Galerkin system~(\ref{galerkin}) of
dimension $mn=5440$.
The spectral abscissa (largest real part of eigenvalues in matrix pencil
$\lambda \hat{E} - \hat{A}$)
becomes $-0.0049$, which confirms the asymptotic stability of the
linear dynamical system. 
We solve the IVP with the input~(\ref{input-periodic}).
The numerical solution yields approximations of the expected value as well as
the standard deviation with respect to the QoI.
Figure~\ref{fig:statistics} illustrates these statistics.

\begin{figure}
\begin{center}
\includegraphics[width=7cm]{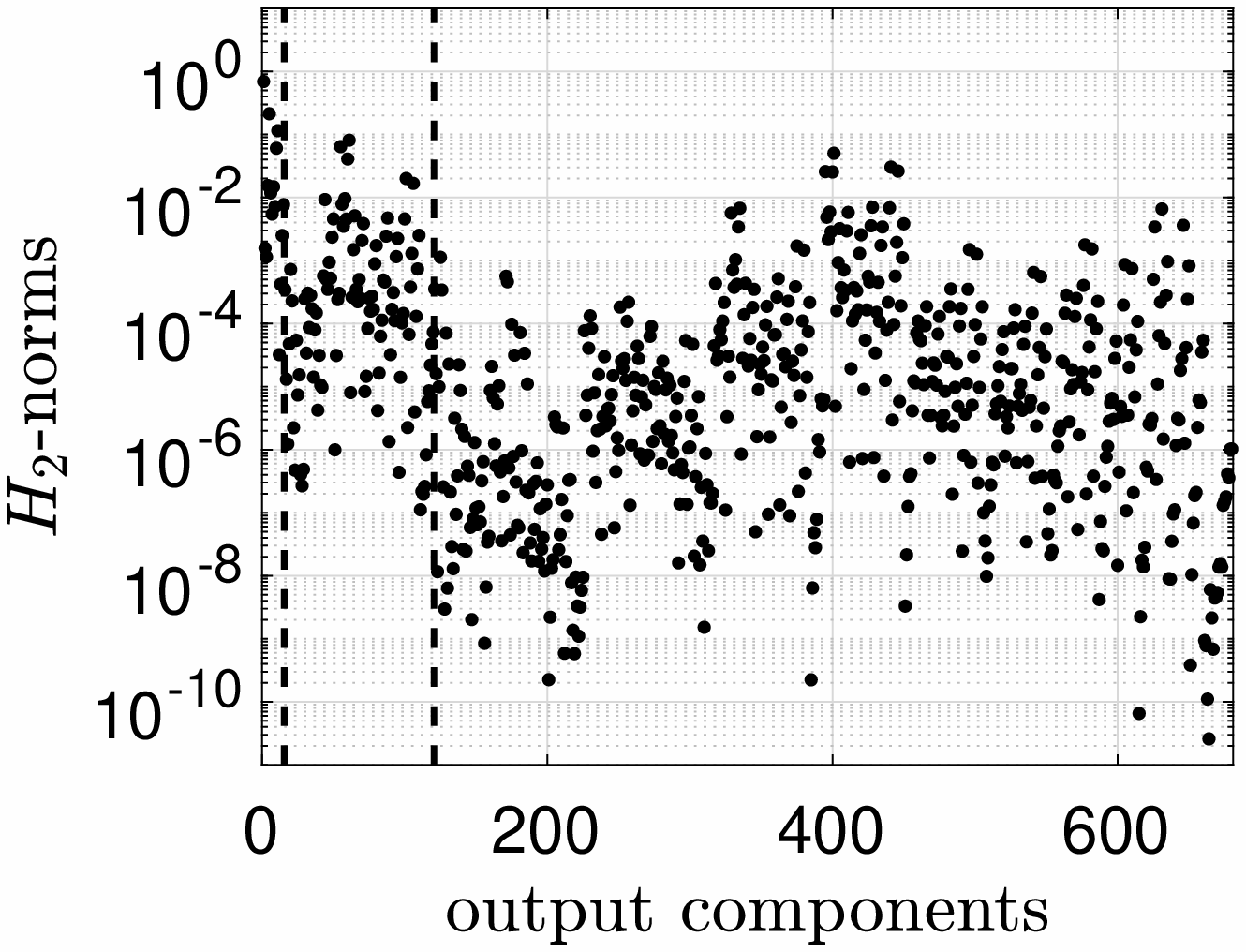}
\hspace{10mm}
\includegraphics[width=7cm]{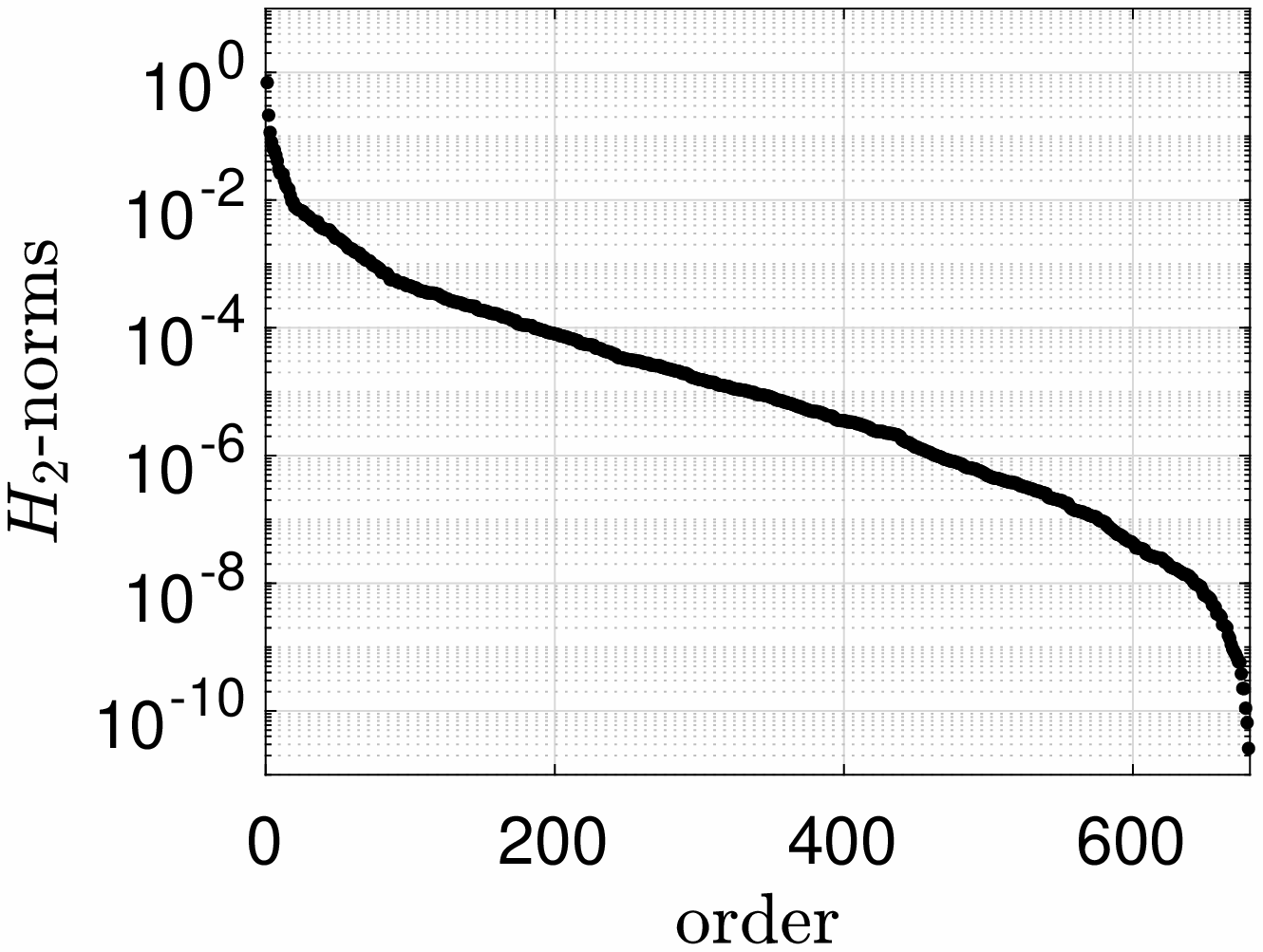}
\end{center}
\caption{$\mathscr{H}_2$-norms of output components in stochastic Galerkin
  system:
  for each component (left),
  dashed lines separate the coefficients for polynomial degree zero/one,
  two and three,
  and in descending order (right).}
\label{fig:h2norms}
\end{figure}

\begin{figure}
\begin{center}
\includegraphics[width=7cm]{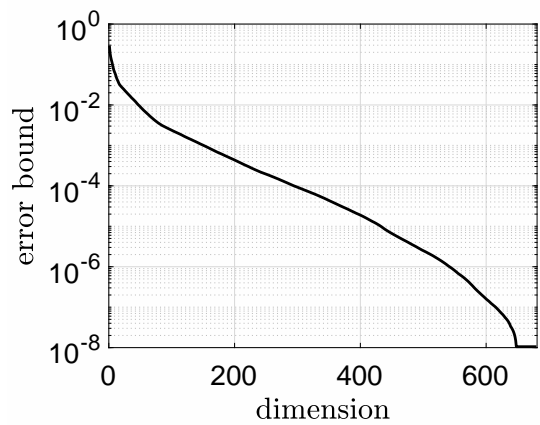}
\end{center}
\caption{Optimal error bounds~(\ref{errorbound}) from $\mathscr{H}_2$-norms
  for input with unit norm given different dimensions of subspaces
  (number of basis polynomials). }
\label{fig:errorbound}
\end{figure}

The quadrature technique produces
approximations of the $\mathscr{H}_2$-norms~(\ref{hardy-norm})
as described in Section~\ref{sec:quadrature}.
We choose $\omega_{\min}=10^{-2}$ and $\omega_{\max}=10^2$.
In the interval~$[\omega_{\rm \min},\omega_{\rm \max}]$,
we use a logarithmically spaced grid with $\nu = 200$ points.
Broader intervals $[\omega_{\min},\omega_{\max}]$ or higher numbers of
grid points do not considerably change the results. 
The computed norms are shown in Figure~\ref{fig:h2norms}.
We observe that the norms exhibit different orders of magnitudes and
a rapid decay.
The $\mathscr{H}_2$-norms allow for the computation of the
error bounds in the right-hand side of~(\ref{errorbound})
using the optimal index sets~(\ref{optimal-index-set})
for increasing dimensions.
Figure~\ref{fig:errorbound} depicts the error estimates
for a normalized input, i.e., $\| u \|_{\ltwotime} = 1$.
Even though the error bounds are not sharp, they decrease exponentially.

%%% Figure: Hardy-norms %%%%%%%%%%%%%%%%%%%%%%%%%%%%%%%%%%%%%%%%%%%%%%%%%%%%%
%%%%%%%%%%%%%%%%%%%%%%%%%%%%%%%%%%%%%%%%%%%%%%%%%%%%%%%%%%%%%%%%%%%%%%%%%%%%%

%%% Figure: Hardy-norms %%%%%%%%%%%%%%%%%%%%%%%%%%%%%%%%%%%%%%%%%%%%%%%%%%%%%
%%%%%%%%%%%%%%%%%%%%%%%%%%%%%%%%%%%%%%%%%%%%%%%%%%%%%%%%%%%%%%%%%%%%%%%%%%%%%

Furthermore, we generate $k=500$ samples of the QoI by solving IVPs
of the linear dynamical system~(\ref{ode})
with input~(\ref{input-periodic}) in the interval~$[0,T]$.
We determine low-dimensional representations~(\ref{low-dim-appr})
in equidistant time points $t_j= j \frac{T}{r}$ for $j=1,\ldots,r$
using $r = 1000$.
The presented results are arithmetic mean values of the quantities
in these time points, which can be seen as approximations of
an integral mean.

\begin{figure}
\begin{center} 
\includegraphics[width=7cm]{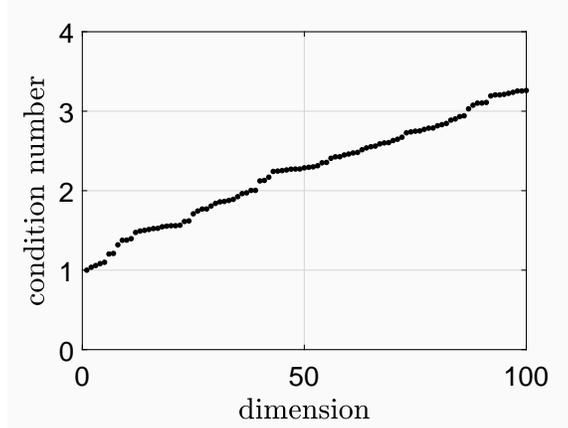}
\end{center}
\caption{Condition numbers of matrices in linear least squares problems~(\ref{least-squares}) for different dimensions.}
\label{fig:condition}
\end{figure}

\begin{figure}
\begin{center}
\includegraphics[width=7cm]{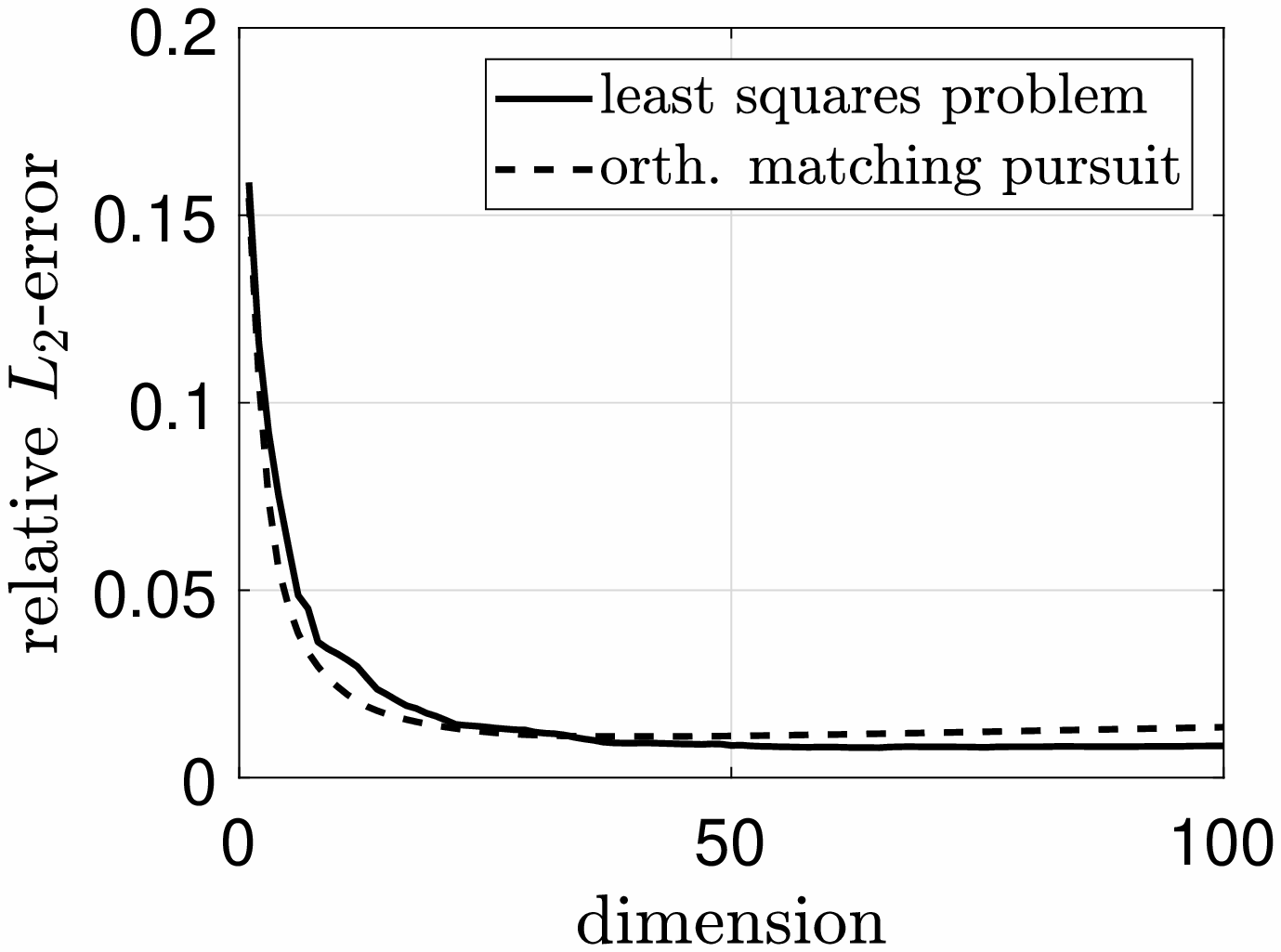}
\hspace{10mm}
\includegraphics[width=7cm]{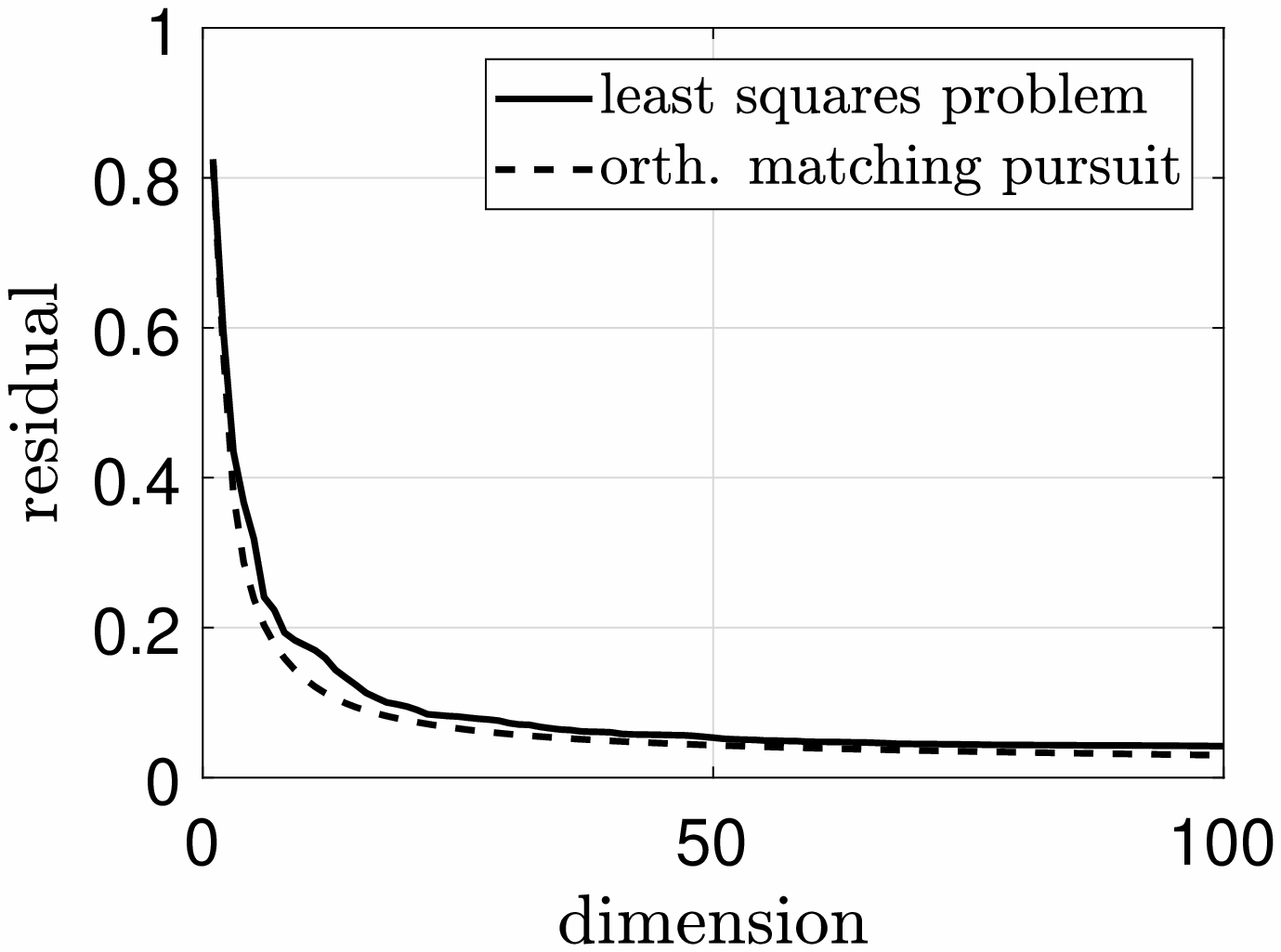}
\end{center}
\caption{Relative $\mathscr{L}^2$-errors~(\ref{l2error}) (left) and
  residuals~(\ref{residual}) (right) of the low-dimensional representations
  obtained by least squares problem and OMP approach.}
\label{fig:error}
\end{figure}

We compute $q$-sparse approximations~(\ref{low-dim-appr})
for $q=1,2,\ldots,100$.
On the one hand, the least squares problems from Section~\ref{sec:lsp}
are solved, where the frequency domain analysis yields the
sequence of bases.
The subspace of dimension~$q$ is time-invariant in this approach.
Figure~\ref{fig:condition} shows the condition numbers
(with respect to the spectral norm) of the matrices
in the least squares problems~(\ref{least-squares}).
We observe that all problems are well-conditioned. 
On the other hand, the minimization problem including OMP
from Section~\ref{sec:l0min} is applied,
where the basis selection takes place separately for each time point.
Concerning the $\mathscr{L}^2$-error, the transient solution of the
Galerkin-projected system~(\ref{galerkin}) serves as reference solution
in the $r$~time points.
Figure~\ref{fig:error} demonstrates both 
$\mathscr{L}^2$-errors~(\ref{l2error}) and residuals~(\ref{residual}).
The residuals decrease monotone in both methods,
which is guaranteed point-wise in time by construction.
The $\mathscr{L}^2$-errors of the least squares problem also
reduce monotone for increasing dimensions of the subspaces.
In the OMP technique, this error decreases rapidly first and then
increases slightly.
The reason is that the underlying minimization problem~(\ref{eq:omp})
with $\varepsilon = 0$ tends to an interpolation of the samples
for higher dimensions of the subspaces.

%%% Figure: Condition numbers %%%%%%%%%%%%%%%%%%%%%%%%%%%%%%%%%%%%%%%%%%%%%%%
%%%%%%%%%%%%%%%%%%%%%%%%%%%%%%%%%%%%%%%%%%%%%%%%%%%%%%%%%%%%%%%%%%%%%%%%%%%%%

%%% Figure: Errors %%%%%%%%%%%%%%%%%%%%%%%%%%%%%%%%%%%%%%%%%%%%%%%%%%%%%%%%%%
%%%%%%%%%%%%%%%%%%%%%%%%%%%%%%%%%%%%%%%%%%%%%%%%%%%%%%%%%%%%%%%%%%%%%%%%%%%%%

\begin{figure}
\begin{center}
\includegraphics[width=7cm]{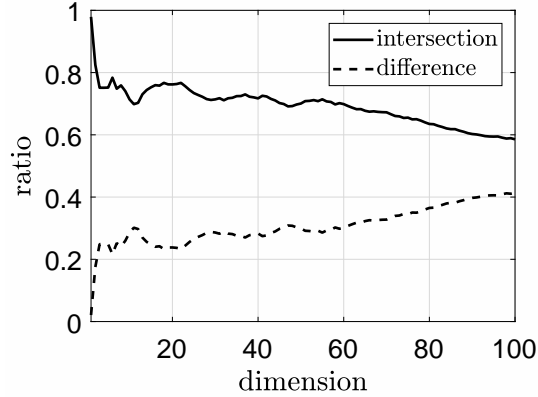}
\end{center}
\caption{Ratio of basis polynomials in the intersection of index sets
  from least squares problem in Section~\ref{sec:lsp} and
  sparse minimization problem in Section~\ref{sec:l0min}.}
\label{fig:ratio-joint}
\end{figure}

For fixed dimension~$q$, we obtain index sets $\mathcal{I}_q^{\rm lsp}$
and $\mathcal{I}_q^{\rm OMP}$ from the two techniques,
which are associated to subsets of basis polynomials.
It holds that $| \mathcal{I}_q^{\rm lsp} | = | \mathcal{I}_q^{\rm OMP} | = q$.
We investigate the agreement of the two sets.
The ratio of basis functions in the intersection reads as
$$ \theta(q) :=
\frac{| \mathcal{I}_q^{\rm lsp} \cap \mathcal{I}_q^{\rm OMP} |}{q} . $$
A ratio $\theta(q)=0$ implies that the two sets are
disjoint, whereas a ratio $\theta(q)=1$ represents identical sets.
The larger the ratio the more the two basis selections agree.
The value $1-\theta(q)$ is the ratio for the differences
$\mathcal{I}_q^{\rm lsp} \backslash \mathcal{I}_q^{\rm OMP}$ as well as
$\mathcal{I}_q^{\rm OMP} \backslash \mathcal{I}_q^{\rm lsp}$.
Figure~\ref{fig:ratio-joint} depicts these ratios,
which represent the mean values in time again. 
We observe that about 60\% or more of the basis functions are
within both index sets for all dimensions.
However, the relative $\mathscr{L}^2$-error, see Figure~\ref{fig:error},
is nearly the same in both methods.
The minimum errors are around $0.01$, which can be seen as a
good low-dimen\-sional approximation.
We conclude that basis polynomials in the intersection of both methods
are the most important contributors to the approximation.

%%% Figure: Ratio of Intersection %%%%%%%%%%%%%%%%%%%%%%%%%%%%%%%%%%%%%%%%%%%
%%%%%%%%%%%%%%%%%%%%%%%%%%%%%%%%%%%%%%%%%%%%%%%%%%%%%%%%%%%%%%%%%%%%%%%%%%%%%

\begin{figure}
\begin{center}
\includegraphics[width=7cm]{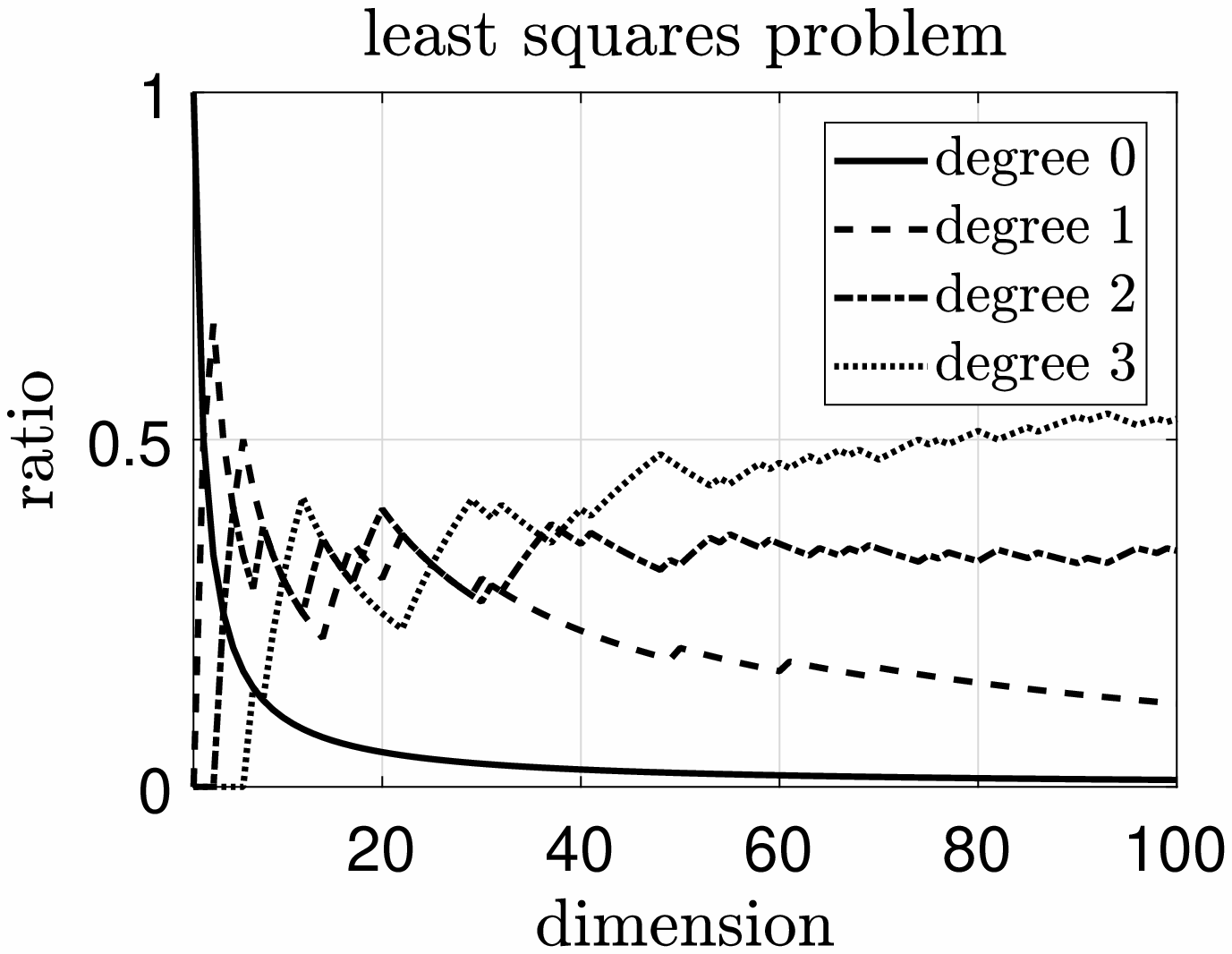}
\hspace{10mm}
\includegraphics[width=7cm]{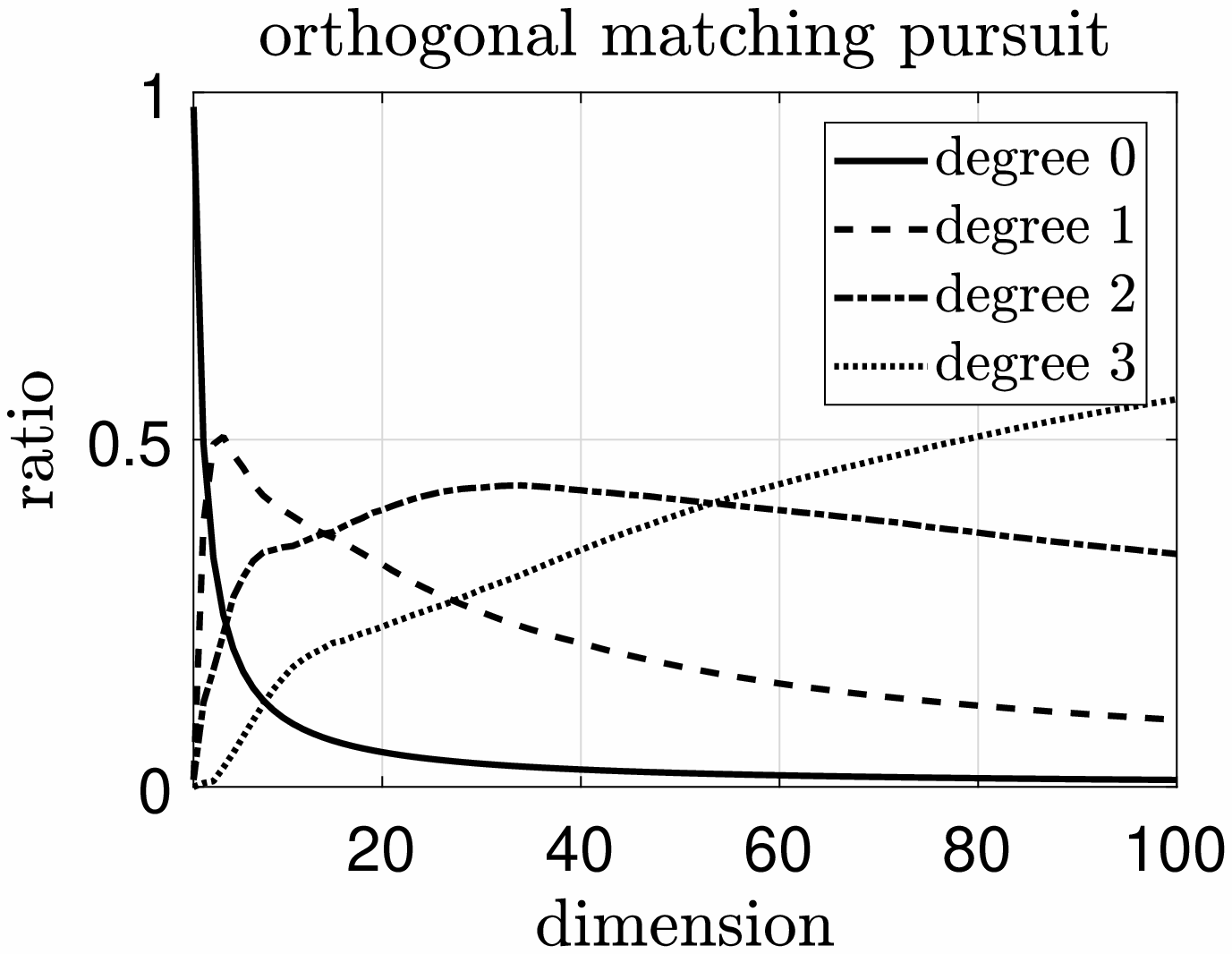}
\end{center}
\caption{Ratios of the polynomials of different degrees in the subspaces
  identified by least squares problem from Section~\ref{sec:lsp} (left) and
  sparse minimization problem from Section~\ref{sec:l0min} (right).}
\label{fig:ratio-degrees}
\end{figure}

In addition, we examine the degrees of the polynomials within
the index sets for the two methods separately.
If $q$~basis polynomials have been chosen in a method, then
let $q_j$ be the number of included polynomials of degree~$j$
($0 \le q_j \le q$). 
We obtain the ratios $0 \le \frac{q_j}{q} \le 1$ for $j=0,1,2,3$.
The sum of the four ratios is equal to one.
Figure~\ref{fig:ratio-degrees} shows that these ratios are similar
in both approaches. 
The polynomial of degree zero is included in nearly all computed bases.
The two techniques are able to identify the dominant basis polynomials
among all degrees.

%%% Figure: Ratio of degrees %%%%%%%%%%%%%%%%%%%%%%%%%%%%%%%%%%%%%%%%%%%%%%%%
%%%%%%%%%%%%%%%%%%%%%%%%%%%%%%%%%%%%%%%%%%%%%%%%%%%%%%%%%%%%%%%%%%%%%%%%%%%%%

%%%%%%%%%%%%%%%%%%%%%%%%%%%%%%%%%%%%%%%%%%%%%%%%%%%%%%%%%%%%%%%%%%%%%%%%%%%%%
%%%                          Conclusions                                  %%%
%%%%%%%%%%%%%%%%%%%%%%%%%%%%%%%%%%%%%%%%%%%%%%%%%%%%%%%%%%%%%%%%%%%%%%%%%%%%%

\section{Conclusions}
We compared two numerical techniques for a basis selection to
achieve low-dimensional approximations in polynomial chaos expansions.
In the first approach 
the basis is identified by an analysis in the frequency domain
followed by a least squares problem for coefficients
in the time domain.
The second approach solves a sparse minimization problem
using information in the time domain. 

We tested the methods on a linear dynamical system modeling
a mechanical configuration.
The numerical results demonstrate that both techniques identify
low-dimensional approximations of the same quality, i.e.,
with a sufficiently small relative $\mathscr{L}^2$-error.
Although the identification of the basis polynomials uses different
approaches, most of the selected basis functions (60\% or more)
coincide in both methods at all times.
Moreover, the number of basis polynomials for each degree separately 
is nearly the same. 

\section*{Acknowledgements}
Sandia National Laboratories is a multimission laboratory managed and
operated by National Technology and Engineering Solutions of Sandia, LLC., a
wholly owned subsidiary of Honeywell International, Inc., for the
U.S. Department of Energy's National Nuclear Security Administration
under contract DE-NA-0003525.
The views expressed in the article do not necessarily represent the views of
the U.S. Department of Energy or the United States Government.
John Jakeman's work was supported by DARPA EQUiPS.

\noindent The authors are indebted to Prof.~Dr.~Akil~Narayan
(University of Utah) for helpful discussions.

% \clearpage

%%%%%%%%%%%%%%%%%%%%%%%%%%%%%%%%%%%%%%%%%%%%%%%%%%%%%%%%%%%%%%%%%%%%%%%%%%%%%
%%%                           References                                  %%%
%%%%%%%%%%%%%%%%%%%%%%%%%%%%%%%%%%%%%%%%%%%%%%%%%%%%%%%%%%%%%%%%%%%%%%%%%%%%%

\end{document}